\documentclass[10pt,reqno]{amsart}
\usepackage{amsmath,amssymb,graphicx}
\usepackage[usenames,dvipsnames]{color}

\numberwithin{equation}{section}
\numberwithin{figure}{section}
\numberwithin{table}{section}
\allowdisplaybreaks[4]

\definecolor{c20}{rgb}{0.,0,0.}
\definecolor{c30}{rgb}{0,0,0}
\definecolor{c40}{rgb}{0,.0,0}
\definecolor{c50}{rgb}{0,0,0}

\def\rrE#1{\textcolor{c30}{#1}}
\def\luAAAA{ \rrE{(\log (u))^4}}
\def\luAA{ \rrE{(\log (u))^2}}

\def\rE#1{\textcolor{c30}{#1}}
\def\rrE#1{\textcolor{c30}{#1}}
\def\eeE#1{\textcolor{c30}{#1}}

\def\pD#1{\textcolor{c40}{#1}}
\def\rrE#1{#1}

\def\rE#1{#1}
\def\pD#1{#1}

\def\cK#1{\textcolor{c40}{#1}}

\def\cK#1{\textcolor{green}{#1}}
\def\cK#1{#1}

\newcommand{\PP}{\mathbb{P}}
\newcommand{\mean}[1]{\mathbb{E}\left( #1\right)}
\newtheorem{theo}{Theorem}[section]
\newtheorem{sat}[theo]{Proposition}
\newtheorem{de}[theo]{Definition}
\newtheorem{lem}[theo]{Lemma}
\newtheorem{exxa}[theo]{Example}

\newtheorem{korr}[theo]{Corollary}

\newtheorem{remarks}[theo]{Remarks}

\newcommand{\nelem}[1]{{Lemma \ref{#1}}}

\newcommand{\netheo}[1]{{Theorem \ref{#1}}}

\newcommand{\prooftheo}[1]{ \textsc{Proof of Theorem} \ref{#1} }

\newcommand{\prooflem}[1]{\textsc{Proof of Lemma} \ref{#1}}

\def\kal#1{{\mathcal{ #1}}}

\newcommand{\E}[1]{\mathbb{E}\left(#1\right)}

\newcommand{\pk}[1]{\mbox{\rm$\mathbb{P}$} \left(#1\right) }

\newcommand{\ldot}{,\ldots,}

\newcommand{\BQN}{\begin{eqnarray}}
\newcommand{\EQN}{\end{eqnarray}}
\newcommand{\BQNY}{\begin{eqnarray*}}
\newcommand{\EQNY}{\end{eqnarray*}}

\newcommand{\BS}{\begin{sat}}
\newcommand{\ES}{\end{sat}}
\newcommand{\BT}{\begin{theo}}
\newcommand{\ET}{\end{theo}}
\newcommand{\BK}{\begin{korr}}
\newcommand{\EK}{\end{korr}}

\newcommand{\BD}{\begin{de}}
\newcommand{\ED}{\end{de}}

\newcommand{\BIT}{\begin{itemize}}
\newcommand{\EIT}{\end{itemize}}
\newcommand{\BDI}{\begin{description}}
\newcommand{\EDI}{\end{description}}

\newcommand{\BRM}{\begin{remarks}}
\newcommand{\ERM}{\end{remarks}}

\newcommand{\QED}{\hfill $\Box$}

\newcommand{\IF}{\infty}
\newcommand{\BTH}{\begin{theo}}
\newcommand{\ETH}{\end{theo}}
\newcommand{\BPR}{\begin{sat}}
\newcommand{\EPR}{\end{sat}}

\newcommand{\BEX}{\begin{exxa}}
\newcommand{\EEX}{\end{exxa}}

\newcommand{\BC}{\begin{cases}}
\newcommand{\EC}{\end{cases}}
\newcommand{\COM}[1]{}
\newcommand{\BL}{\begin{lem}}
\newcommand{\EL}{\end{lem}}

\begin{document}

\begin{center}
\thispagestyle{empty}

{\large \bf Tail Asymptotics of Random Sum and Maximum of Log-Normal Risks}

       \vskip 0.4 cm

         \centerline{\large
         Enkelejd Hashorva\footnote{University of Lausanne, UNIL-Dorigny 1015 Lausanne, Switzerland} and Dominik Kortschak\footnote{
Universit\'e de Lyon, F-69622, Lyon, France; Universit\'e Lyon 1, Laboratoire SAF, EA 2429, Institut de Science Financi\`ere et d'Assurances, 50 Avenue Tony Garnier, F-69007 Lyon, France}
}

 \vskip 0.4 cm

\end{center}

{\bf Abstract:} In this paper we derive the asymptotic behaviour of the survival function of both random sum and random maximum of log-normal risks. As for the case of finite sum and maximum investigated in  Asmussen and Rojas-Nandaypa (2008) also for the more general setup of random sums and random maximum the principle  of a single big jump holds. We investigate both the log-normal sequences and some related dependence structures motivated by stationary Gaussian sequences.


\section{Introduction}
Let $Y_i,i\ge 1$ be positive random variables \rE{(rv's)} which model claim sizes of an insurance portfolio for a given observation period.
Denote by $N$ the total number of \rrE{claims} reported during the observation period, thus $N$ is a discrete rv, which we assume to be independent of claim sizes $Y_i,i\ge 1$. The classical \rE{risk} model $S_N= \sum_{i=1}^N Y_i$ for the total loss amount assumes that $Y_i$'s are independent and identically distributed (iid) rv's.
If the assumption of independence of claim sizes is dropped, \pD{one faces the problem how to choose a meaningful dependence structure.  Further this dependence structure should be tractable from a theoretical point of view. For example Constantinescu et al.\ (2011)  consider \rE{a} model \rE{where}  the survival copula of claim sizes is assumed to be Archimedean. Such a model
 has the  interpretation  that for some positive rv $V$
 and  iid unit exponential rv's $E_i,i\ge 1$  independent of $V$, \rE{then} $Y_i=V E_i,i\ge 1$ form a dependent sequence of claim sizes
 derived by randomly scaling of iid claim sizes $E_i,i\ge 1$.}\\
\pD{In this paper we use dependent Gaussian sequences and related dependence structures \rE{to model claim sizes}.}
 Specifically, if $X_i,i\ge 1$ \rE{are} dependent Gaussian rv's with $N(0,1)$ distribution, then $Y_i= e^{X_i},i\ge 1$ is the corresponding sequence of dependent log-normal rv's that can be used for modeling claim sizes. For instance, if $X_i,i\ge 1$ is a centered stationary Gaussian sequence of $N(0,1)$ components and constant correlation $\rho= \E{X_1X_i}\in (0,1),i>1$, then $Y_i=e^{X_i}$ is a sequence of dependent log-normal rv's.
 Since we have (see e.g., Berman (1992))
\BQN
X_i= \rho Z_0+ \sqrt{1- \rho^2} Z_i, \label{seq}
\EQN
with $Z_i,i\ge 0$ iid $N(0,1)$ rv's, then $Y_i= e^{\rho Z_0} e^{\sqrt{1- \rho^2}Z_i},i\ge 1$. For such $Y_i$'s,
by Asmussen and Rojas-Nandaypa (2008)
\BQN\label{NN}
\pk{S_n > u} \sim n \pk{X_1> \log u}, \quad u\to \IF
\EQN
holds for any $n\ge 2$, where $\sim$ stands for asymptotic equivalence of two functions when the argument tends to infinity. In view of Asmussen et al.\ (2011) (see also Hashorva (2013)) $S_n$ is asymptotically tail equivalent with the maximum $Y_{n:n}= \max_{1 \le  i \le n} Y_i$, i.e., $\pk{S_n> u} \sim \pk{Y_{n:n} >u}$ as $u\to \IF$.

Our analysis in this paper is concerned with the probability of observing large values for the random sum $S_N$, thus we shall investigate $\pk{S_N> u}$ when $u$ is large. Additionally, we shall consider also the tail asymptotics of the maximum claim $Y_{N:N}$ among
the claim sizes $Y_1 \ldot Y_N$; we set $Y_{0:0}=0$ if $N=0$. \rrE{For the case that $N$ is non-random see for recent results on max-sum equivalence  Jiang et al.\ (2014) and the references therein.}

For our investigations of the tail behaviours of $S_N$ and $Y_{N:N}$ we shall follow two objectives:\\
A) We shall exploit the tractable
dependence structure \rE{implied} by \eqref{seq} choosing general  $Z_i$'s such that $e^{Z_i}$ has survival function similar to that of a log-normal
rv; \\
B) We consider  a log-normal dependence structure induced by a general Gaussian sequence $X_i,i\ge 1$ where $X_i,X_n$ can have a correlation $\rho_{in}$ which is allowed to converge to 1 as $n\to \IF$.

For both cases of dependent $Y_i$'s we show that the principle of a single big jump (see Foss et al.\ (2013) for details in iid setup)
holds if for the discrete rv $N$
we require that
\BQN \label{conN}
 \E{(1+ \delta)^N} < \IF
\EQN
is valid for some $\delta>0$; a large class of discrete rv's satisfies condition \eqref{conN}.\\
Brief organisation of the rest of the paper: We present our main results in Section 2
followed by the proofs in Section 3.

\section{Main Results}
We consider first $X_i$'s which \rrE{are in general not Gaussian}. So for a given fixed $\rho \in [0,1)$ let $Z_i,i\ge 0$ be independent rv's which define $X_i$'s via the dependence structure \eqref{seq}. We shall assume that 
\BQN\label{14b}
\pk{e^{Z_0} > u}\sim \mathcal{L}(u) \Psi(\log(u)), \quad u\to \IF,
\EQN
 with  $\Psi$ the survival function of an $N(0,1)$ rv and $\mathcal{L}(\cdot)$ a regularly varying function at $\IF$ with index $\beta\in \mathbb{R}$, see Bingham et al. (1987) or Mikosch (2009) for details on regularly \rrE{varying} functions.  Clearly, \eqref{14b} is satisfied if $Z_0$ is an $N(0,1)$ rv. Considering $Z_0$ as a base risk, we shall further assume that with $c_i\in [0,\IF)$ uniformly in $i$
\BQN\label{14}
\pk{Z_i > u}\sim c_i\pk{{Z_0} > u}, \quad u\to \IF.
\EQN
For such models the claim sizes
$Y_i=e^{X_i},i\ge 1$  have marginal distributions which are in general neither log-normal nor with tails which are proportional \rrE{to} those of log-normal rv's.

We state next our first result \rE{for $Y_{N:N}$ the maximal claim size among $Y_1=e^{X_1}  \ldot Y_N=e^{X_N}$
 and the random sum $S_N=\sum_{i=1}^N Y_i$; we set $Y_{0:0}=0$ and $S_0:=0$.}

\BT\label{th0}  Let $N$ be an integer-valued rv satisfying $\E{(1+\delta)^N} < \IF$ for some $\delta>0$. Let $X_i,i\ge 1$ be a sequence of rv's given by \eqref{seq}
with $Z_i,i\ge 0$ iid rv's and $\rho\in [0,1)$ some given constant.  Suppose that \eqref{14b} and \eqref{14}  hold
 with  $\max_{i \ge 1} c_i < \IF$.  \rE{If further $N$ is independent of $X_i,i\ge 1$, then}
\BQN\label{NN2}
\pk{S_N > u} \sim \pk{Y_{N:N} > u} \sim \E{\sum_{i=1}^N c_i}  \frac{\kal{L}(u^{\rho^2}) \kal{L}(u^{1-\rho^2}) }{ \sqrt{ 2 \pi} \log u} \exp\Bigl( - \frac{(\log u)^2}{2} \Bigr) , \quad u \to \IF.
\EQN
\ET

{\bf Remarks:} a) Clearly, if $Y=e^{Z}$ with $Z$ an $N(0,1)$ rv (thus $Y$ is a log-normal rv with $LN(0,1)$ distribution), then \eqref{14b} holds with $\kal{L}(u)=1, u>0$. \\
b) If $\kal{L}(\cdot)$ in \netheo{th0} is constant, then the tail asymptotic behaviour of $S_N$ and $Y_{N:N}$ is not influenced by the value of the dependence parameter $\rho$, and hence as expected the principle of a
single big jump  holds. \rrE{However}, for non-constant $\kal{L}(\cdot)$ the dependence  parameter $\rho$ plays a crucial  role in the tail asymptotics derived in \eqref{NN2}. The reason for this is that by Lemma \ref{L1}
\begin{equation}
 \pk{Y_i>u} \sim c_i \frac{\kal{L}(u^{\rho^2}) \kal{L}(u^{1-\rho^2}) }{ \sqrt{ 2 \pi} \log u} \exp\Bigl( - \frac{(\log u)^2}{2} \Bigr),
 \quad u\to \IF.
 \label{ff}
\end{equation}
Hence also in this case the principle of a single big jump applies.\\
c) In the proof of Theorem \ref{th0} we can show $S_N\stackrel{d}{=}e^{\rho Z_0}e^{\sqrt{1-\rho^2} Z^\star}$  for some $Z^*$ independent of $Z_0$ and then we apply Lemma \ref{L1}. Here we  want to mention that after proving
\eqref{ff} we can also apply Proposition 2.2 of Foss and Richards (2010) to determine the asymptotic of $\pk{S_n> u}$ as $u\to \IF$.
If we condition on $Z_0$ and set  $\overline F(x) =\PP(Y_1>x)$,
$B_i(x)=\{x:e^{\rho Z_0}\le x^\gamma\}$ for some $\gamma  \in (\rho,1)$ and define $\pD{h(x)=x^{\xi}}$ with \[1-\frac 12 \left(\frac{1-\gamma}{\sqrt{1-\rho^2}}\right)<\xi^2<1,
\]
then it is  straightforward to show that the conditions of Proposition 2.2 of Foss and Richards (2010) are met.

\bigskip

Our second result is for log-normal rv's where we remove the assumptions of equi-correlations.
Specifically, we consider for each $n$ claim sizes
$Y_{1,n}=e^{X_{1,n}},\ldots,Y_{n,n}= e^{X_{n,n}}, $ where $(X_{1,n}, \ldots, X_{n,n})$ is a normal random vector with mean zero and covariance matrix $\Sigma^{(n)}$ which is a correlation matrix with entries $\sigma^{(n)}_{i,j}$. We shall assume that $\rho^n_{i,j}:=\sigma^{(n)}_{i,j}$ is bounded by some sequence $\rho_n$ and some $\rho\in (0,1)$, i.e.,
\BQN\label{rhon}
\rho^n_{i,j}\le\max(\rho_n,\rho), \quad \rE{n\ge 1}
\EQN
for all  $i\not =j$. \rrE{Further, we suppose that the} sequence  $\rho_n,n\ge 1$ \rrE{satisfies} for some $c^*>8$ and some $\eta> 0$
\BQN\label{cnu}
 \rho_{n(u)} \le 1-\frac{c^* \log(\log(u))}{\log(u)}, \quad \text{with }
 n(u)= \left \lfloor(1+ \eta ) \frac{ \luAA }{2\log(1+\delta)}\right\rfloor .
\EQN
If for instance all \rE{${\rho_{i,j}^{n}}$} are bounded, then clearly condition \eqref{cnu} is valid; it holds also if for some $c$ large enough $\rho_n\le 1-c \log(n)/\sqrt{n}$.\\
We present next our final result.
\begin{theo}\label{thm}
Let $Y_{1,n} \ldot Y_{n,n}, n\ge 1$ be claim sizes as above being further independent of some integer-valued rv $N$ which satisfies
\eqref{conN} for some $\delta>0$. If further \eqref{rhon} holds with $\rho_n$ satisfying \eqref{cnu}, then 
\BQN
 \PP\left( \max_{1\le i \le N}Y_{i,N}>u\right)\sim  \PP(S_N>u)\sim
\rE{\frac{ \mean{N} }{ \sqrt{ 2 \pi} \log u} \exp\Bigl( - \frac{(\log u)^2}{2} \Bigr), \quad u\to \IF}.
\EQN
\end{theo}

{\bf Remarks:}
a) Our second result in \netheo{thm} shows that the principle of a single big jump  still holds even if we allow \rrE{for a more general dependence structure}.\\
b) Kortschak (2012) derives second order asymptotic results for subexponential risks. Similar ideas as therein are utilised to derive second order asymptotic results for the aggregation of log-normal random vectors in Kortschak and Hashorva (\rrE{2013,2014}). In the setup of randomly weighted sums it is also possible to derive such results.

\def\OF{\overline{F}}

\section{Proofs}
We give next  two lemmas needed in the proofs below. \rrE{The first lemma is of some interest on its own, 
in particular it implies Lemma 2.3 in Farkas and Hashorva (2013) (see also Lemma 8.6 in Piterbarg (1996)).}

\def\td{ \eta} 

\BL \label{L1} Let  $\kal{L}_i\rE{(\cdot)}, i=1,2$ be
some regularly varying functions at infinity with index $\beta_i$.
If   $Z_1,Z_2$ are two independent rv such that  $\pk{e^{Z_i}> u} \sim \kal{L}_i(u) \Psi(\log(u)), i=1,2$, then
for any $\sigma_1,\sigma_2$ two positive constants
\BQN
\pk{e^{\sigma_1 Z_1+\sigma_2 Z_2} >u} \sim  \sigma^2 e^{\frac{\sigma_1^2 \sigma_2^2}{2 \sigma^2} (\beta_1 - \beta_2)^2 } \kal{L}_1(u^{\gamma})\kal{L}_2(u^{1-\gamma}) \Psi( (\log u)/\sigma)
\label{38}
 \EQN
holds as $u\to \IF$, where $\gamma=\sigma_1^2/(\sigma_1^2+\sigma_2^2)$ and $\sigma=\sqrt{\sigma_1^2+ \sigma_2^2}$.
\EL
\def\sY{\sigma_*}
\def\sYd{\sigma_*^2}

\prooflem{L1}
Choose an $\alpha>0$ such that \[
                              \frac{\sigma_1^2}{\sigma_2^2}<\frac{1+\alpha}{1-\alpha}.
                             \]
Then for  any $a>0$ we have
\BQN
\frac{ \pk{e^{\sigma_1 Z_1 +\sigma_2 Z_2 }> u,  e^{\sigma_2 Z_2}  \le a}}{\pk{e^{\sigma_1 Z_1+ \sigma_2 Z_2 }> u,  e^{\sigma_2 Z_2}> a}} & \le &
\frac{ \pk{e^ { \sigma_1 Z_1 } > u/a}} {\pk{e^{\sigma_1 Z_1} > u^\alpha} \pk{{e^{\sigma_2Z_2}}> u^{1-\alpha}}}\notag\\
&\sim & \frac{ \kal{L}_1( (u/a)^{1/\sigma_1})\Psi(\frac{1}{\sigma_1}\log (u/a))}{ \kal{L}_1( u^{\alpha/\sigma_1}) \kal{L}_2( u^{(1-\alpha)/\sigma_1}) \Psi(\frac{\alpha}{\sigma_1}\log (u)) \Psi(\frac{1-\alpha}{\sigma_2}\log (u)) }\notag\\
& \to & 0, \quad u\to \IF,\label{eq:temp1}
\EQN
with $\Psi$ the survival function of an $N(0,1)$ rv.  With the same argument we get that
 for  any $a>0$ we have
\BQNY
\frac{ \pk{e^{\sigma_1 Z_1 +\sigma_2 Z_2 }> u,  e^{\sigma_1 Z_1}  \le a}}{\pk{e^{\sigma_1 Z_1+ \sigma_2 Z_2 }> u,  e^{\sigma_1 Z_1}> a}} \to 0, \quad u\to \IF,\label{eq:temp1C}
\EQNY
and hence 
\begin{align}
 \pk{e^{\sigma_1 Z_1 +\sigma_2 Z_2 }> u}=&\pk{e^{\sigma_1 Z_1 +\sigma_2 Z_2 }> u,  e^{\sigma_1 Z_1}  > a,  e^{\sigma_2 Z_2}  > a}\notag\\&+\pk{e^{\sigma_1 Z_1 +\sigma_2 Z_2 }> u,  e^{\sigma_1 Z_1}  \le a}+\pk{e^{\sigma_1 Z_1 +\sigma_2 Z_2 }> u,  e^{\sigma_2 Z_2}  \le a}\notag
\\ 
\sim&\pk{e^{\sigma_1 Z_1 +\sigma_2 Z_2 }> u,  e^{\sigma_1 Z_1}  > a,  e^{\sigma_2 Z_2}  > a}, \quad u\to \IF.\label{eq:temp1b}
\end{align}
\COM{For $\xi_u=u^{1-\gamma }/\xi,\xi>0$ we have $u/\xi_u= \xi u^{-\gamma}$
\BQNY
\pk{ e^{\sigma_1 Z_1 }> u/ \xi_u }  \pk{e^{\sigma_2 Z_2}> \xi_u}&= &
\pk{ e^{\sigma_1 Z_1 }> \xi u^\gamma }  \pk{e^{\sigma_2 Z_2}> u^{1- \gamma}/\xi}\\
&=& \pk{  Z_1 > (\gamma \ln u + \ln \xi)/\sigma_1} \pk{Z_2 > ((1- \gamma) \ln u - \ln \xi)/\sigma_2}.
\EQNY
}
In view of \eqref{eq:temp1b} we have 
\BQNY
\pk{e^{\sigma_1 Z_1+ \sigma_2 Z_2} > u}
&\sim& \pk{e^{\sigma_1 Z_1+ \sigma_2 Z_2} > u,e^{\sigma_1 Z_1}>\xi,e^{\sigma_2 Z_2}>\xi}, \quad u\to \IF.
\EQNY
Assume next without loss of generalty that $\sigma_1 \ge \sigma_2$. 
If  $H$ denotes the distribution of $e^{\sigma_1 Z_1}$, then for any $\xi>0$ with $u>2\xi$ 
\BQNY
\pk{e^{\sigma_1 Z_1+ \sigma_2 Z_2} > u,e^{\sigma_1 Z_1}>\xi,e^{\sigma_2 Z_2}>\xi}
&=&\pk{e^{\sigma_1 Z_1+ \sigma_2 Z_2} > u,u/\xi \ge e^{\sigma_1 Z_1}>\xi,e^{\sigma_2 Z_2}>\xi}\\&&+\ \pk{e^{\sigma_1 Z_1+ \sigma_2 Z_2} > u,e^{\sigma_1 Z_1}>u/\xi,e^{\sigma_2 Z_2}>\xi}\\
&=&\pk{e^{\sigma_1 Z_1+ \sigma_2 Z_2} > u,u/\xi \ge e^{\sigma_1 Z_1}>\xi}+\pk{e^{\sigma_1 Z_1}>u/\xi,e^{\sigma_2 Z_2}>\xi}\\
&=&  \int_{\xi}^{u/\xi} \pk{ e^{\sigma_2 Z_2}> u/s} \ d H(s)+  \pk{ e^{\sigma_1 Z_1 }> u/ \xi ,e^{\sigma_2 Z_2}> \xi}.
\EQNY
For all  $u $ and $\xi$ large enough  
$$\int_{\xi}^{u/\xi} \pk{ e^{\sigma_2 Z_2}> u/s} \ d H(s)\ge 
\frac{1}{2}\pk{e^{\sigma_2 Z_2}> u/\xi}\ge  \pk{ e^{\sigma_1 Z_1 }> u/ \xi ,e^{\sigma_2 Z_2}> \xi}$$
implying as $u\to \IF$
\BQNY
\int_{\xi}^{u/\xi} \pk{ e^{\sigma_2 Z_2}> u/s} \ d H(s)+  \pk{ e^{\sigma_1 Z_1 }> u/ \xi ,e^{\sigma_2 Z_2}> \xi}
&\sim&   \int_{\xi}^{u/\xi} \pk{ e^{\sigma_2 Z_2}> u/s} \ d H(s).
\EQNY
Further, since again the constant $\xi$ can be chosen arbitrary large we get for $\gamma=\sigma_1^2/(\sigma_1^2+\sigma_2^2)$
\BQNY
\lefteqn{\int_{\xi}^{u/\xi} \pk{ e^{\sigma_2 Z_2}> u/s} \ d H(s)} \\
&\sim&
   \int_{\xi}^{u/\xi} \frac{\sigma_2^2\kal{L}_2(u/s)} {\sqrt{2 \pi  \sigma_2^2}\log  (u/s)} \exp\Biggl( - \frac{(\log  (u/s))^2}{2  \sigma_2^2}\Biggr)\ d H(s)\\
&=&  \frac{\sigma_2^2\kal{L}_2\left(u^{1-\gamma}\right)} {\sqrt{2 \pi  \sigma_2^2}\log  (u^{1- \gamma})} \int_{\xi u^{-\gamma}}^{\frac 1\xi u^{1-\gamma}} \frac{\kal{L}_2\left(\frac{u^{1-\gamma}} s\right) } {\kal{L}_2\left(u^{1-\gamma}\right) }
\frac{\log \left(u^{1-\gamma}\right) } {\log \left(\frac{u^{1-\gamma}} s\right) }
 \exp\Biggl( - \frac{\left(\log  \left(\frac{u^{1-\gamma}} s \right)\right)^2}{2  \sigma_2^2}\Biggr) \ d H( u^{\gamma}s)\\
&=&  \frac{(\sigma_1^2+\sigma_2^2)\kal{L}_2\left(u^{1-\gamma}\right)} {\sqrt{2 \pi  \sigma_2^2}\log  (u)} \int_{\xi u^{-\gamma}}^{\frac 1\xi u^{1-\gamma}}
q(u,  \gamma,s) \exp\Biggl( - \frac{\left(\log  \left(\frac{u^{1-\gamma}} s \right)\right)^2}{2  \sigma_2^2}\Biggr) \ d H( u^{\gamma}s),
\EQNY
with $q(u,\gamma,s)=\frac{\kal{L}_2\left(\frac{u^{1-\gamma}} s\right) } {\kal{L}_2\left(u^{1-\gamma}\right) }
\frac{\log \left(u^{1-\gamma}\right) } {\log \left(\frac{u^{1-\gamma}} s\right) }.$  For some $c>0$, \rE{by the uniform convergence theorem for regularly varying functions (see Theorem A3.2 in Embrechts et al.\ (1997))}  we get uniformly in  $1/c<s<c$
\[
\lim_{u\to\infty} q(u,\gamma,s)
=s^{-\beta_2}.
\]
\rrE{Further note that in the light of Potter's bound \rE{(see Bingham et al.\ (1987))}} for every $\epsilon>0$ and $A>1$ we can find  a positive constant \eeE{$\xi$} such that for  all $\xi u^{-\gamma} <s< \frac 1\xi u^{1-\gamma}$
\[
\frac 1 A s^{-\beta_2} \min(s^\epsilon,s^{-\epsilon}) \le
q(u,\gamma,s) \le A s^{-\beta_2} \max(s^\epsilon,s^{-\epsilon}).
\]
Consequently, for different values of $0 <a<b$ (that might depend on $u$) and $\beta$ 
we want to find the asymptotics of
\begin{align*}
& \int_{a}^{b}
s^\beta
\exp\Biggl( - \frac{\left(\log  \left(\frac{u^{1-\gamma}} s \right)\right)^2}{2  \sigma_2^2}\Biggr) \ d H( u^{\gamma}s)\\&= -s^\beta
\exp\Biggl( - \frac{\left(\log  \left(\frac{u^{1-\gamma}} s \right)\right)^2}{2  \sigma_2^2}\Biggr)\PP\left(e^{\sigma_1Z_1}>u^{\gamma}s \right) \Bigg|_{s=a}^b \\&
 \quad
+\int_a^b  s^{\beta-1} \left(\beta + \frac{\log\left(\frac{u^{1-\gamma}} s \right)}{\sigma_2^2} \right)
\exp\Biggl( - \frac{\left(\log  \left(\frac{u^{1-\gamma}} s \right)\right)^2}{2  \sigma_2^2}\Biggr)
\PP\left(e^{\sigma_1Z_1}>u^{\gamma}s \right) \ ds.
\end{align*}
Since we can choose $\xi$ arbitrary large we can replace $\PP\left(e^{\sigma_1Z_1}>u^{\gamma}s \right)$ by its asymptotic form and hence we can use the approximation (set $\sY:= \sigma_1 \sigma_2/\sqrt{\sigma_1^2+ \sigma_2^2}$)
\begin{align*}
&\exp\Biggl( - \frac{\left(\log  \left(\frac{u^{1-\gamma}} s \right)\right)^2}{2  \sigma_2^2}\Biggr)
\PP\left(e^{\sigma_1Z_1}>u^{\gamma}s \right)\\&
\approx \sigma_1^2\frac{\kal{L}_1(u^\gamma s)} {\sqrt{2 \pi  \sigma_1^2}\log  (u^\gamma s)} \exp\Biggl( - \frac{\left(\log  \left(\frac{u^{1-\gamma}} s \right)\right)^2}{2  \sigma_2^2}-\frac{\left(\log  \left(u^{\gamma}s \right)\right)^2}{2  \sigma_1^2}\Biggr)\\
&= \sigma_1^2\frac{\kal{L}_1(u^\gamma s)} {\sqrt{2 \pi  \sigma_1^2}\log  (u^\gamma s)}
\exp\Biggl( - \frac{\left(\sigma_1^2(1-\gamma)^2+\sigma_2^2\gamma^2\right)  (\log (u))^2 + 2\left(\sigma_1^2(\gamma-1)+\sigma_2^2\gamma\right) \log (u)\log (s)+(\sigma_1^2+\sigma_2^2)(\log (s))^2  }{2  \sigma_1^2\sigma_2^2}\Biggr)\\
&=\sigma_1^2
\frac{\kal{L}_1(u^\gamma s)} {\sqrt{2 \pi  \sigma_1^2}\log  (u^\gamma s)}
\exp\Biggl( - \frac{(\log (u))^2}{2(\sigma_1^2+\sigma_2^2)}\Biggr) \exp\Biggl( -\frac{(\log (s))^2  }{2  \sYd}\Biggr).
\end{align*}
Since $ \sigma_1^2(\gamma -1)+ \sigma_2^2 \gamma=0$, using again Potter's bounds \rE{(see Bingham et al.\ (1987))} and the fact that $\kal{L}_1(\cdot)$ is regularly varying at infinity,
the above derivations imply
\begin{align*}
& \PP(e^{\sigma_1 Z_1+\sigma_2 Z_2}>u)\\ &\sim   \frac{\sigma_1^2 (\sigma_1^2 +\sigma_2^2)  \kal{L}_1(u^{\gamma}) \kal{L}_2(u^{1-\gamma})} {\sigma_2^2\sqrt{2\pi \sigma_2^2\sigma_1^2} \log (u)}\frac{1-\gamma}{\gamma\sqrt{2\pi} }
\exp\Biggl( - \frac{(\log (u))^2}{2(\sigma_1^2+\sigma_2^2)}\Biggr)
\int_0^\infty s^{\beta_1-\beta_2-1} \exp\Biggl( -\frac{(\log (s))^2  }{2 \sYd }\Biggr) \ ds  \\
&= \frac{ \sqrt{\sigma_1^2+\sigma_2^2} \kal{L}_1(u^{\gamma}) \kal{L}_2(u^{1-\gamma})} {\sqrt{2 \pi}  \log (u)}
\exp\Biggl( - \frac{(\log (u))^2}{2(\sigma_1^2+\sigma_2^2)}\Biggr)
\int_0^\infty \frac{1} {\sqrt{2\pi  \sY^2}}s^{\beta_1-\beta_2-1} \exp\Biggl( - \frac{(\log (s))^2  }{2  \sYd}\Biggr) \ ds\\
&= \sqrt{\sigma_1^2+\sigma_2^2} e^{\frac{\sYd }{2} (\beta_1-\beta_2)^2 }\frac{ \kal{L}_1(u^{\gamma}) \kal{L}_2(u^{1-\gamma})} {\sqrt{2 \pi}  \log (u)}
\exp\Biggl( - \frac{(\log (u))^2}{2(\sigma_1^2+\sigma_2^2)}\Biggr),
\end{align*}
hence the proof is complete. \QED

\begin{lem}\label{lem:asymptotictwo} Assume that $n\le n(u)$ with $n(u)$ \rE{defined in \eqref{cnu}}  and \rrE{set} $\epsilon(u)=4\log(\log(u))/\log(u)$. If $Y_1$ is an $LN(0,1)$ rv and $X_{i,n}, i\le n$ are as in \netheo{thm},  then  as $u\to \IF$
\[
\PP(Y_1>u-n u^{1-\epsilon(u)}) \sim \PP(Y_1>u)
\]
and for $i\not=j$
\[
\PP(Y_{i,n}>u^{1-\epsilon(u)} ,Y_{j,n}>u^{1-\epsilon(u)}) =o(\PP(Y_1>u)).
\]
\end{lem}
\prooflem{lem:asymptotictwo}  By the assumptions on $n$ and $n(u)$ as $u\to \IF$ we have
\begin{align*}
\PP(Y_1>u)\le \PP(Y_1>u-n u^{1-\epsilon(u)}) &\le \PP(Y_1>u-n(u) u^{1-\epsilon(u)})\\
&=\PP\left(Y_1>u-\frac u{  \luAAAA} (1+ \eta ) \frac{ \luAA }{2\log(1+\delta)}\right)\\
&=\PP\left(Y_1>u-\frac{(1+ \eta )}{2\log(1+\delta)} \frac u{ \luAA }  \right)\\&\sim\PP(Y_1>u).
\end{align*}
\rE{Next, denote by  $f$ the probability density function of $Y_1$}. Let further
$\rrE{W_1} $ and $\rrE{W_2} $ be two independent $N(0,1)$ rv's,  and write  \rE{$\rho_*$} for the correlation between \rE{$\log Y_{i,n}$ and $\log Y_{j,n}$}.
We may write for $u>0$
\BQNY
\PP(Y_{i,n}>u^{1-\epsilon(u)} ,Y_{j,n}>u^{1-\epsilon(u)})
&=&
\PP( e^{\rrE{W_1} }> u^{1-\epsilon(u)}, e^{\rE{\rho_*}\rrE{W_1} +\sqrt{1-\rE{\rho_*}^2} \rrE{W_2} } > u^{1-\epsilon(u)})\\
&=& \PP\left(e^{\rrE{W_1} }>\frac{u}{ \luAAAA},e^{ \rE{\rho_*}\rrE{W_1} } e^{\sqrt{1-\rE{\rho_*}^2} \rrE{W_2} } >\frac{u}{ \luAAAA  } \right)\\
&\le &\PP\left(\frac u{  \luAAAA  }<e^{\rrE{W_1} }<2u,e^{ \rE{\rho_*} \rrE{W_1} } e^{\sqrt{1-\rE{\rho_*}^2} \rrE{W_2} } >\frac{u}{ \luAAAA  } \right)+\PP(e^{\rrE{W_1} }>2u)\\
&=&
 \int_{\frac{u}{ \luAAAA  }}^{2u}
  \PP\left(e^{\rrE{W_2} }>\left(\frac{u}{ \luAAAA  x^{\rE{\rho_*}}}\right)^{1/\sqrt{1-\rE{\rho_*}^2}} \right)
 f(x) d x +\PP(e^{\rrE{W_1} }>2u)\\
&\le &\int_{\frac{u}{ \luAAAA  }}^{2u}
  \PP\left(e^{\rrE{W_2} }>\left(\frac{u^{1-\rho_*}}{ \luAAAA  2^{\rho_*}}\right)^{1/\sqrt{1-\rE{\rho_*}^2}} \right)
 f(x) d x+\PP(e^{\rrE{W_1} }>2u)\\
&\le &\PP\left(Y_1>\frac{u^{\sqrt{\frac{1-\rho_*}{1+\rho_*}}}}{2^{\rho_*} \luAAAA  }\right)\PP\left(Y_1>\frac{u}{ \luAAAA  } \right)+\PP(e^{\rrE{W_1} }>2u)\\
&=&o(\PP(Y_1>u)), \quad u\to \IF
\EQNY
since 
\BQNY
\left(1+\frac{1-\rho_*}{1+\rho_*}\right)\log(u)&=&\frac{2}{1+\rho_*}\log(u)\\
&\ge& \frac{2}{1+\rho_{n(u)}}\log(u)\\
& \ge& \frac{2}{2-\frac{c^* \log(\log(u))} {\log(u)}} \log(u)\\
&=&\log(u) +\frac{2}{2-\frac{c^* \log(\log(u))} {\log(u)}} c^* \log(\log(u))\\
& \sim &\log(u) +c^* \log(\log(u)).
\EQNY
Consequently, the assumption  $c^*>8$ entails 
\begin{align*}
& 2\log\left(\PP\left(Y_1>\frac{u^{\sqrt{\frac{1-\rho_*}{1+\rho_*}}}}{2^{\rho_*} \luAAAA  }\right)\PP\left(Y_1>\frac{u}{ \luAAAA  } \right) \right)\\& \sim \log\left(\frac{u^{\sqrt{\frac{1-\rho_*}{1+\rho_*}}}}{2^{\rho_*} \luAAAA  }\right)^2+
\log\left(\frac{u}{\luAAAA  }\right)^2\\
&\sim \left(1+\frac{1-\rho_*}{1+\rho_*}\right)\log(u)^2-8\log(u) \log(\log(u))-8\sqrt{\frac{1-\rho_*}{1+\rho_*}}\log(u) \log(\log(u))\\
&\lesssim \log(u)^2+(c^*-8) \log(\log(u))
\end{align*}
establishing the proof. \QED

\prooftheo{th0}
For any $u>0$ we have
\BQNY
 \pk{S_N> u}&=& \pk{ e^{\rho Z_0} \sum_{i=1}^N e^{ \sqrt{1- \rho ^2} Z_i} > u}\\
 &=:& \pk{ e^{\rho Z_0} W_N  > u}.
\EQNY
Since $e^{ \sqrt{1- \rho^2} Z_i},i\ge 1$ are subexponential risks, then  along the lines of the proof of
Theorem 3.37 in Foss et al.\ (2013) (see also for similar result Theorem 1.3.9 in Embrechts et al.\ (1997))
$$\pk{W_N> u} \sim \Theta\pk{e^{ \sqrt{1- \rho^2} Z^*}> u}, \quad \Theta:= \E{\sum_{i=1}^N c_i} $$
as $u\to \IF$, with $Z^*$ an independent copy of $Z_0$. It can be easily checked that  $Z_0$ and $\log(W_N)/(1-\rho^2)$ fulfill the conditions of  \nelem{L1}, hence
the asymptotic of  $\PP(S_N>u)$  follows. Similarly, 
\BQNY
 Y_{N:N}&=& \max_{1 \le i \le N} e^{ \rho Z_0+ \sqrt{1- \rho^2} Z_i}
 =  e^{\rho Z_0 }  \max_{1 \le i\le N} e^{ \sqrt{1 -\rho^2} Z_i}=:e^{\rho Z_0} W_N^*.
 \EQNY
Since we have
$$ \pk{W_N^*> u} \sim \pk{W_N>u} \sim \Theta \pk{ \exp( \sqrt{1 - \rho^2} Z^*)> u}, \quad u \to \IF$$
the proof follows by applying once again \nelem{L1}.  \QED

\prooftheo{thm} Denote next $Y_1$ an $LN(0,1)$ rv and let $\mathcal{I}_{\{\cdot\}}$ denote the indicator function.
Since for all fixed $n\ge 1$ we get by interchanging limit and finite sum that
\BQNY
\PP(S_N>u)&=&\PP(S_N>u,N\le n)+\PP(S_N>u,N> n)\\
&\sim& \mean{N \mathcal I_{\{ N\le n\}}} \PP(Y_1>u)+\PP(S_N>u,N> n)
\EQNY
we can assume w.l.o.g. that $\rho_{i,j}^n\le \rho_n$.  From \eqref{conN} it follows that there exist $C_1,C_2>0$ such that
\begin{align*}
  p_n:=\PP(N=n)\le C_1 (1+\delta)^{-n} \quad\text{and}\quad \PP(N>n)\le C_2(1+\delta)^{-n}.
\end{align*}
By the independence of $N$ and the claim sizes
\[
 \PP(S_N>u)=\sum_{n=1}^\infty p_n\PP(S_n>u)
\]
and for $n(u)$ defined in \eqref{cnu}
\BQNY
 \sum_{n=n(u)}^\infty p_n\PP(S_n>u)&\le& \PP(N>n(u)) \\
 &\le &C_2(1+\delta)^{-n(u)} \\
 &\le &C_2\exp\left(-\frac{1+ \eta }2 \luAA \right)\\
 &=&o(\PP(Y_1>u)).
\EQNY
Since
\[\PP(S_n>u)\ge n\PP(Y_1>u)-\sum_{i\not=j} \PP(Y_i>u,Y_j>u)\]
and by Lemma \ref{lem:asymptotictwo}
$$\PP(Y_i>u,Y_j>u)=o(\PP(Y_1>u)), \quad u\to \IF$$
 it follows that
\BQNY
  \sum_{n=0} ^{n(u)} p_n\PP(S_n>u)&\ge& \PP(Y_1>u) \left(\sum_{n=0} ^{n(u)} n p_n- o(1)\sum_{n=0} ^{n(u)} n^2 p_n \right)\\
  &\sim &\mean{N} \PP(Y_1>u), \quad u\to \IF.
\EQNY
So we are left with finding an asymptotic upper bound. For $n\le n(u)$ we use the \rrE{following} decomposition (c.f. \pD{Asmussen  and Rojas-Nandaypa (2008)})
\begin{align*}
 \PP(S_n>u)=\sum_{i=1}^n \PP\left(S_n>u,Y_{i,n}\ge Y_{j,n},\max_{j\not=i} Y_{j,n}>u^{1-\epsilon(u)}\right) +\PP\left(S_n>u,Y_{i,n}\ge Y_{j,n},\max_{j\not=i} Y_{j,n}\le u^{1-\epsilon(u)}\right),
\end{align*}
where $\epsilon(u)=4\log(\log(u))/\log(u)$. By Lemma \ref{lem:asymptotictwo} we have
\begin{align*}
 \sum_{i=1}^n \PP\left(S_n>u,Y_{i,n}\ge Y_{j,n},\max_{j\not=i} Y_{j,n}>u^{1-\epsilon(u)}\right)
&\le \sum_{i=1}^n\sum_{j\not=i} \PP\left(S_n>u,Y_{i,n}\ge Y_{j,n}, Y_{j,n}>u^{1-\epsilon(u)}\right)\\
&\le \sum_{i=1}^n\sum_{j\not=i} \PP\left(Y_{i,n}>u^{1-\epsilon(u)},Y_{j,n}>u^{1-\epsilon(u)}\right)\\
&=n(n-1) o(\PP(Y_1>u)).
\end{align*}
Further 
\begin{align*}
 \PP\left(S_n>u,Y_{i,n}\ge Y_{j,n},\max_{j\not=i} Y_{j,n}\le u^{1-\epsilon(u)}\right)&\le\PP\left(Y_{i,n}>u-\sum_{i=1}^n Y_{j,n},\max_{j\not=i} Y_{j,n}\le u^{1-\epsilon(u)}\right)\\
&\le \PP\left(Y_{i,n}>u-n u^{1-\epsilon(u)},\max_{j\not=i} Y_{j,n}\le u^{1-\epsilon(u)}\right)\\
&\le\PP\left(Y_{i,n}>u-n u^{1-\epsilon(u)}\right) \\
& \sim \PP(Y_1>u)
\end{align*}
as $u \to \IF$, hence  the proof for the tail asymptotics of $S_N$ follows \eeE{by} applying \eqref{lem:asymptotictwo}. Since for any $u>0$
\[ n\PP(Y_1>u)-\sum_{i\not=j} \PP(Y_i>u,Y_j>u)\le \PP\left( \max_{1\le i \le n}Y_{i,n}>u\right)\le  \PP(S_n>u)
\]
\rrE{the tail asymptotics of $\max_{1\le i \le N}Y_{i,N}$ can be easily established, \eeE{and} thus the proof is complete.}
 \QED

\textbf{Acknowledgments.} We would like to that the referees of the paper for several suggestions which improved our manuscript.
E. Hashorva kindly acknowledges partially support  by the Swiss
National Science Foundation Grant 200021-140633/1 and RARE -318984 (an FP7 Marie Curie IRSES Fellowship).

\end{document}